\documentclass[12pt]{amsart}%
\usepackage{amsmath}
\usepackage{amsfonts}
\usepackage{amssymb}
\usepackage{graphicx}%
\newtheorem{theorem}{Theorem}[section]
\theoremstyle{plain}

\newtheorem{condition}[theorem]{Condition}

\newtheorem{corollary}[theorem]{Corollary}

\newtheorem{definition}[theorem]{Definition}
\newtheorem{example}[theorem]{Example}

\newtheorem{lemma}[theorem]{Lemma}

\newtheorem{problem}[theorem]{Problem}
\newtheorem{proposition}[theorem]{Proposition}
\newtheorem{remark}[theorem]{Remark}

\numberwithin{equation}{section}
\begin{document}
\title[Bounded Lipschitz Metric]{Bounded--Lipschitz Distances on the State Space of a C*--algebra}
\author{Fr\'{e}d\'{e}ric Latr\'{e}moli\`{e}re}
\address{Department of Mathematics, University of Toronto}
\email{frederic@math.toronto.edu}
\urladdr{http://www.math.toronto.edu/frederic}

\begin{abstract}
Metric noncommutative geometry, initiated by Alain Connes, has known some
great recent developments under the impulsion of Rieffel and the introduction
of the category of compact quantum metric spaces topologized thanks to the
quantum Rieffel-Gromov-Hausdorff distance. In this paper, we undertake the
first step to generalize such results and constructions to locally compact
quantum metric spaces. Our present work shows how to generalize the
construction of the bounded-Lipschitz metric on the state space of a
C*-algebra which need not be unital, such that the topology induced by this
distance on the state space is the weak* topology. In doing so we obtain some
results on a state space picture of the strict topology of a C*-algebra.

\end{abstract}
\maketitle

\section{Introduction}

In noncommutative geometry, as suggested in \cite{Connes89},\cite[Ch.
6]{Connes}, the natural way to specify a (quantum) metric on a C*-algebra $A$
is to choose a densely-defined seminorm $L$ on $A$, which should be viewed as
a generalization of a Lipschitz seminorm. Given such a seminorm $L$ and any
$\alpha>0$, one can define a norm $L_{\alpha}=\max\left\{  L,\alpha
^{-1}\left\Vert .\right\Vert \right\}  $ whose dual seminorm on $A^{\prime}$
then induces a distance $d_{L,\alpha,1}$ on the state space $S(A)$ of $A$. The
prototype for this situation is given by a locally compact separable metric
space $\left(  X,\rho\right)  $ and by taking $L$ to be the usual Lipschitz
seminorm on $A=C_{0}(X)$. The natural embedding of $X$ into the space
$S(C_{0}(X))$ of Radon probability measures over $X$, which maps any $x\in X$
to the Dirac point measure $\delta_{x}$, is also an isometry from any balls in
$X$ of radius at most $\alpha$ into $\left(  S(C_{0}(X)),d_{L_{\alpha}%
}\right)  $ for any $\alpha>0$. Therefore, the family of norms $\left(
L_{\alpha}\right)  _{\alpha>0}$ encodes all the metric information of $X$ just
as the C*-norm encodes all the topological information of $X$. In this
context, the distances $d_{L,\alpha,1}$ are known as the bounded-Lipschitz
distances and have the fundamental property that any one of them induces the
weak* topology on $S(C_{0}(X))$ \cite[Theorem 11.3.3 p. 395]{Dudley}. A
natural question, therefore, is to find out what condition on $L$ one needs to
ensure that $d_{L,\alpha,1}$ metrizes the weak* topology on $S(A)$ for any
(non-Abelian) C*-algebra $A$.

When $A$ is unital, this problem was answered by Rieffel in \cite[Theorems 1.8
and 1.9]{Rieffel98a}: the distance $d_{L,\alpha,1}$ induces the weak* topology
on $S(A)$ if and only if the set $\left\{  a\in A:L(a)\leq1,\left\Vert
a\right\Vert \leq\alpha\right\}  $ is norm precompact in $A$. This generalizes
the prototypical situation where $X$ is assumed to be compact, in which case
the norm-precompactness of any set of $1$-Lipschitz uniformly bounded
functions follows from Arz\'{e}la-Ascoli. Moreover, when $X $ is compact, it
is of finite diameter $r$, so that as long as $L(1)=0$, the bounded-Lipschitz
distance $d_{L,r,1}$ is also the distance induced by the dual seminorm of $L$
and is known in this context as the Kantorovich distance. It is in this
framework that \cite[Theorems 1.8 and 1.9]{Rieffel98a} is written. Following
on \cite{Rieffel98a}, it is possible to develop a theory of Gromov-Hausdorff
convergence \cite{Rieffel00} for C*-algebras with such seminorms, leading to
many interesting new approximations of C*-algebras, as in \cite{Rieffel00},
\cite{Rieffel01},\cite{Latremoliere05}.

We undertake in this paper the first step into a generalization of the theory
of compact quantum metric spaces, as developed in \cite{Rieffel98a}%
,\cite{Rieffel99} and \cite{Rieffel00}, to non-unital C*-algebras. As in the
unital case, this first step consists in characterizing those seminorms on
(nonunital) C*-algebras whose associated bounded-Lipschitz distances induces
the weak* topology on $S(A)$. We thus answer the problem:

\begin{problem}
\label{Pb}Let $A$ be a separable C*-algebra, with or without a unit. Let $L$
be a seminorm on a dense subset $\operatorname*{dom}L$ of the self-adjoint
part $A^{\operatorname*{sa}}$ of $A$. We extend $L$ to $A^{\operatorname*{sa}%
}$ by setting $L(a)=\infty$ whenever $a\not \in \operatorname*{dom}L$. We
define the distance $d_{L}=d_{L,1,1}$ on $S(A)$ by setting $d_{L}\left(
\varphi,\psi\right)  =\sup\left\{  \left\vert \varphi(a)-\psi(a)\right\vert
:L(a)\leq1\text{, }\left\Vert a\right\Vert \leq1\right\}  $ for all
$\varphi,\psi\in S(A)$.When does the $d_{L}$-topology agree with the weak*
topology on $S(A)$?
\end{problem}

The choice of the constant $\alpha=1$ in Problem (\ref{Pb}) can be made
without loss of generality, as shown in Proposition (\ref{metricseq}), since
we are only asking a topological, i.e. local question.

The choice of bounded-Lipschitz distances, rather than the Kantorovich
distance $\kappa_{L}=d_{L,\infty,1}$, as the natural framework when working
with nonunital C*-algebras is justified by the following two simple
observations about $\kappa_{L}$ when $A$ is Abelian nonunital and $L$ is the
Lipschitz seminorm for some distance on the Gelfan'd spectrum $X$ of $A$.
First of all, $\kappa_{L}$ is valued in $[0,\operatorname*{diam}X]$ and thus
can take the value $\infty$ if $X$ has infinite radius:

\begin{example}
\label{real_line0}Let $X=\mathbb{R}$ with its usual metric and denote by $L$
the usual Lipschitz seminorm on $C_{0}\left(  \mathbb{R}\right)  $. Set
$\varphi=\sum_{n=0}^{\infty}2^{-n}\delta_{2^{2n}}\in S(C_{0}(X))$. For any set
$E\subseteq\mathbb{R}$ we define $\chi_{E}$ to be the indicator function of
$E$. Now, for any $n\in\mathbb{N}$ we define $f_{n}(t)=\chi_{\lbrack0,2^{2n}%
[}t+\chi_{\lbrack2^{2n},2^{2n+1}]}\left(  2^{2n+1}-t\right)  $ for all
$t\in\mathbb{R}$. Of course, $L(f_{n})=1$, yet $\varphi(f_{n})\geq2^{n}$.
Since $\delta_{0}(f_{n})=0$ we have $\kappa_{L}(\varphi,\delta_{0})\geq2^{n}$
for any $n\in\mathbb{N}$ so $\kappa_{L}(\varphi,\delta_{0})=\infty$.
\end{example}

Of course, one could still define a topology based on a distance taking the
value $\infty$. However, when $X$ is noncompact, $\kappa_{L}$ does not in
general metrize the weak* topology:

\begin{example}
\label{real_line1}We work with the notations of Example (\ref{real_line0}).
For any $n\in\mathbb{N}\backslash\{0\}$ we define $g_{n}(t)=\chi_{\lbrack
0,n[}t+\chi_{\lbrack n,2n]}(2n-t)$ and $\varphi_{n}=\left(  \left(  1-\frac
{1}{n}\right)  \delta_{0}+\frac{1}{n}\delta_{n}\right)  $. For all
$n\in\mathbb{N}\backslash\{0\}$, we have $\varphi_{n}(g_{n})=1$ and
$L(g_{n})=1$ by construction, so this shows that $\kappa_{L}(\delta
_{0},\varphi_{n})\geq1$ since $\delta_{0}(g_{n})=0$. On the other hand, for
any $f\in C_{0}(\mathbb{R)}$ we have $\varphi_{n}(f)\underset{n\rightarrow
\infty}{\longrightarrow}\delta_{0}(f)$ so $\left(  \varphi_{n}\right)
_{n\in\mathbb{N}}$ weak* converge to $\delta_{0}$.
\end{example}

Therefore, we will focus in this paper on the bounded-Lipschitz distances,
which are known to always metrize the weak* topology for locally compact
spaces, and as we saw include $\kappa_{L}$ when $X$ is compact.

We adopt the point of view that Problem (\ref{Pb}) is of intrinsic interest,
since the state space $S(A)$ is in general a very delicate object (see for
instance \cite{Alfsen71}), especially when $A$ is not unital. Thus, solving
Problem (\ref{Pb}) would provide a mean to bring back problems involving the
weak*$\ $topology $S(A)$ to problems on the C*-algebra $A$ with a well-chosen
seminorm $L$. This approach led to successful developments when $A$ is unital
(\cite{Rieffel00}, \cite{Latremoliere05} to quote a few) and we hope that
further research will show similar results for the nonunital case as well,
based upon our answer to Problem (\ref{Pb}).

A key element in the proof of \cite[Theorems 1.8 and 1.9]{Rieffel98a} is that
$S(A)$ is weak* compact when $A$ is unital. Of course, when $A$ is not unital,
$S(A)$ is never compact for the weak* topology, as its weak* closure $Q(A)$ is
the space of positive linear functionals of norm at most 1 \cite{Pedersen79}.
Moreover, since $S(A)\ $is not weak*\ open in $Q(A)$ (as $Q(A)/S(A)\ $is
dense), it is not even weak*\ locally compact. Therefore we need new
topological insights to solve Problem (\ref{Pb}) in a manner similar to
\cite[Theorems 1.8 and 1.9]{Rieffel98a}.

We propose a solution to Problem (\ref{Pb}) which again provides a topological
condition on $\mathcal{B}_{L}=\left\{  a\in A^{\operatorname*{sa}}%
:L(a)\leq1,\left\Vert a\right\Vert \leq1\right\}  $ which mirrors the topology
of $S(A)$: in the first part of this paper, we show that the $d_{L}$-topology
is the weak* topology on $S(A)\ $if and only if $\mathcal{B}_{L}$ is totally
bounded for the topology \textrm{wu} of uniform weak convergence on weak*
compact subsets of $S(A)$. Since the topology \textrm{wu} is still rather
obscure, we construct in the second part of this paper a metrization for
\textrm{wu }on bounded subsets and compare \textrm{wu} to the strict topology,
which is always stronger than \textrm{wu }and agrees with \textrm{wu} on any
bounded subset $B$ of $A$ if and only if the multiplication is jointly
continuous for \textrm{wu} on $B$, which is usually not the case. We obtain in
this process an alternative description of the strict topology on bounded
subsets using $S(A)$. We then conclude with an easier answer to Problem
(\ref{Pb}) as a corollary.

\bigskip The author wishes to thank Marc Rieffel for his help and support
during this research.

\section{Metrizing the Weak* topology of the state space of a C*-algebra}

In this section, $A$ is a separable C*--algebra with or without a unit, and
whose norm is denoted by $\left\Vert .\right\Vert $. Its state space
$S(A)\ $is endowed with the restriction of the weak*\ topology $\sigma
(A^{\prime},A) $. The self-adjoint part of $A$ is denoted by
$A^{\operatorname*{sa}}$.

\begin{definition}
\label{BL}Let $A$ be a separable C*-algebra and $L$ be a seminorm defined on a
norm-dense subset $\operatorname*{dom}L$ of the set $A^{\operatorname*{sa}}$
of self-adjoint elements in $A$. We extend $L$ to $A^{\operatorname*{sa}}$ by
setting $L(a)=\infty$ for all $a\not \in \operatorname*{dom}L$, and we set:
\[
\mathcal{B}_{L}=\left\{  a\in\operatorname*{dom}L:L(a)\leq1\text{, }\left\Vert
a\right\Vert \leq1\right\}  \text{.}%
\]

\end{definition}

We fix in this section a seminorm $L$ as in Definition (\ref{BL}). The
following easy lemma will prove useful for defining a distance on the state
space $S(A)$ of $A$:

\begin{lemma}
Suppose that, for some $\varphi,\psi\in S(A)$ we have $\varphi(a)=\psi(a)$ for
all $a\in\mathcal{B}_{L}$. Then $\varphi=\psi$.
\end{lemma}

\begin{proof}
Let $a\in\operatorname*{dom}L$. Then $\left(  \max\{1,L(a),\left\Vert
a\right\Vert \right)  ^{-1}a\in\mathcal{B}_{L}$. Hence, $\varphi(a)=\psi(a) $.
Now, $\operatorname*{dom}L$ is norm dense in $A^{\operatorname*{sa}}$, so by
continuity $\varphi=\psi$ on $A^{\operatorname*{sa}}$. Since
$A^{\operatorname*{sa}}$ linearly spans $A$, we have $\varphi=\psi$ on $A$ by linearity.
\end{proof}

\begin{definition}
\label{boundedlipdef}Let $A$ be a separable C*-algebra and $L$ be a densely
defined seminorm on the set $A^{\operatorname*{sa}}$ of self-adjoint elements
of $A $. We define the bounded-Lipschitz distance $d_{L}$ dual to $L$ on the
state space $S(A)$ of $A$ by setting, for all $\varphi,\psi\in S(A)$:
\[
d_{L}(\varphi,\psi)=\sup\left\{  \left\vert \varphi(a)-\psi(a)\right\vert
:a\in\mathcal{B}_{L}\right\}  \text{.}%
\]

\end{definition}

\bigskip We aim at solving Problem (\ref{Pb}). We start with a natural
necessary condition on $L$ for $d_{L}$ to metrize $\sigma(A^{\prime},A)$.
Given a seminorm $p$ on $A$, we shall say a subset $V$ of $A$ is totally
bounded for $p$ when for any $\varepsilon>0$ there exists a finite subset
$F_{\varepsilon}$ of $A$ such that $V\subseteq\{a\in A:p(a)\leq\varepsilon
\}+F_{\varepsilon}$.

\begin{lemma}
\label{QuasiLipNeed}If the $d_{L}$-topology is the weak* topology on
$S(A)\ $then, for all weak* compact subsets $K$ of $S(A)$, the set
$\mathcal{B}_{L}$ is totally bounded for the seminorm $p_{K}$ on $A$ defined
for all $a\in A$ by $p_{K}(a)=\sup_{\varphi\in K}\left\vert \varphi
(a)\right\vert $.
\end{lemma}

\begin{proof}
Assume that $d_{L}$ metrizes the topology $\sigma(A^{\prime},A)$ on $S(A)$.
Let $K$ be a weak* compact subset of $S(A)$.For any $a\in A$ we denote by
$\widehat{a}$ the affine function $\widehat{a}:\varphi\in K\longmapsto
\varphi(a)$. The map $\widehat{a}$ is continuous on $K$ by definition of the
weak* topology. Let $C(K)$ be the Banach space of continuous functions on $K$
endowed with the supremum norm. The map $\xi_{K}:a\longmapsto\widehat{a}\in
C(K)$ is an isometry for $p_{K}$ by definition. Since $p_{K}\leq\left\Vert
.\right\Vert $, the map $\xi_{K}$ is also continuous from $A$ into $C(K)$. We
have, by definition of $d_{L}$, for all $\varphi,\psi\in K$ and for all
$a\in\mathcal{B}_{L}$:
\[
\left\vert \widehat{a}(\varphi)-\widehat{a}(\psi)\right\vert =\left\vert
\varphi(a)-\psi(a)\right\vert \leq d_{L}(\varphi,\psi)\text{.}%
\]
Now, $d_{L}$ metrizes $S(A)$ for the weak* topology, and so in particular the
topology on $K$ obtained from $d_{L}$ is the restriction of $\sigma(A^{\prime
},A)$ on $K$. Therefore, for all $a\in\mathcal{B}_{L}$, the functions
$\widehat{a}$ are equicontinuous over the weak*\ compact $K$. They also are,
by definition, valued in the common compact $[-1,1]$. By Arzela-Ascoli,
$\mathcal{B}_{L}$ is totally bounded for the seminorm $p_{K}$.
\end{proof}

This proof is rather classic, and can be found in \cite[Theorem 1.8]%
{Rieffel98a} when $A$ is a unital\ C*-algebra, in which case the state space
$S(A)$ is itself weak* compact and the functional representation $a\in
A^{\operatorname*{sa}}\longmapsto\widehat{a}\in\operatorname*{Af}\left(
S(A)\right)  $ is a surjective linear order-preserving isometry (see
\cite[sec. 3.10, pp. 69--73]{Pedersen79}) known as the Kadison representation
\cite{Kadisson51}, where $\operatorname*{Af}\left(  S(A)\right)  $ is the
space of continuous affine functions on the convex $S(A)$.

The proof of Lemma (\ref{QuasiLipNeed}) now suggests that we define the
following topology, in search of the proof of a sufficient condition:

\begin{definition}
\label{wu}Let $A$ be a separable C*-algebra. We denote by $\mathfrak{S}$ the
set of all weak* compact subsets of the state space $S(A)$ of $A$. For any
$K\in\mathfrak{S}$ we define the seminorm $p_{K}$ on $A$ by:
\[
p_{K}(a)=\sup_{\varphi\in K}\left\vert \varphi(a)\right\vert \text{\ \ for all
}a\in A\text{.}%
\]
The locally convex topology on $A$ generated by the family $(p_{K}%
)_{K\in\mathfrak{S}}$ of seminorms is called the weakly-uniform topology on
$A$ and is denoted by \textrm{wu}.
\end{definition}

Thus, the topology \textrm{wu} is the topology of uniform weak convergence on
weak* compact subsets of $S(A)$. Note that a subset of $S(A)$ is weak* compact
in $A^{\prime}$ if and only if it is weak* compact in $S(A)$.

We will study the main properties of \textrm{wu} in the next section. This
topology is the natural setting not only for the necessary condition of Lemma
(\ref{QuasiLipNeed}), but also for the sufficient condition in the following
theorem which is the first step toward answering Problem (\ref{Pb}):

\begin{theorem}
\label{quasilip}Let $L$ be a seminorm defined on a norm-dense subset
$\operatorname*{dom}L$ of the set $A^{\operatorname*{sa}}$ of self-adjoint
elements of a separable C*-algebra $A$. Let $\mathcal{B}_{L}=\left\{
a\in\operatorname*{dom}L:L(a)\leq1,\left\Vert a\right\Vert \leq1\right\}  $.
The bounded-Lipschitz distance $d_{L}$ dual to $L$ metrizes the weak* topology
on the state space $S(A)$ of $A$ if, and only if, the set $\mathcal{B}_{L}$ is
totally bounded for the weak-uniform topology \textrm{wu}.
\end{theorem}

\begin{proof}
We established most of the necessary conditions in Lemma (\ref{QuasiLipNeed}).
Let us thus assume that $d_{L}$ metrizes the weak*\ topology on $S(A)$. Let
$V$ be any open neighborhood of $0$ for the $\mathrm{wu}$\textrm{\ }topology.
There exists a set $K\in\mathfrak{S}$ and $\varepsilon>0$ such that $\left\{
a:p_{K}(a)<\varepsilon\right\}  \subseteq V$. By Lemma (\ref{QuasiLipNeed}),
there exists a finite subset $F$ of $A$ such that $\mathcal{B}_{L}%
\subseteq\left\{  a:p_{K}(a)<\varepsilon\right\}  +F\subseteq V+F$. Hence,
$\mathcal{B}_{L}$ is \textrm{wu}--totally bounded.

\qquad Conversely, assume that $\mathcal{B}_{L}$ is $\mathrm{wu}$--totally
bounded. Let $(\varphi_{n})_{n\in\mathbb{N}}$ be a sequence in $S(A)$
converging in $\sigma(A^{\prime},A)$ to a \emph{state} $\varphi_{\infty}$.
Since $K=\{\varphi_{n}:n\in\mathbb{N\cup\{\infty\}}\}\in\mathfrak{S}$, by
assumption $\mathcal{B}_{L}$ is totally bounded for $p_{K}$. Let
$\varepsilon>0$ and let $V$ be the \textrm{wu-}open set $\{a\in A:p_{K}%
(a)<\frac{\varepsilon}{3}\}$. There exists a finite subset $F=\left\{
a_{1},\ldots,a_{\#F}\right\}  $ of $A$ such that $\mathcal{B}_{L}\subseteq V+F
$. Let $N$ be chosen such that for all $n\geq N$, for all $i=0,\ldots,\#F $,
we have $\left\vert \varphi_{n}(a_{i})-\varphi_{\infty}(a_{i})\right\vert
\leq\frac{1}{3}\varepsilon$. Let $a\in\mathcal{B}_{L}$. There exists $a_{i}\in
F$ such that $a-a_{i}\in V$. Therefore, for all $n\geq N$:
\begin{align*}
\left\vert \varphi_{n}(a)-\varphi_{\infty}(a)\right\vert  &  \leq\left\vert
\varphi_{n}(a-a_{i})\right\vert +\left\vert \varphi_{\infty}(a-a_{i}%
)\right\vert +\left\vert \varphi_{n}(a_{i})-\varphi_{\infty}(a_{i})\right\vert
\\
&  \leq\frac{2}{3}\varepsilon+\frac{1}{3}\varepsilon=\varepsilon\text{.}%
\end{align*}

Hence $d_{L}(\varphi_{n},\varphi_{\infty})\leq\varepsilon$. Since $A$ is
separable, the topology $\sigma(A^{\prime},A)$ is metrizable on $S(A)$ by
Lemma\ (\ref{Gdelta}), so this concludes our proof.
\end{proof}

\begin{remark}
\label{quasilip0}Theorem (\ref{quasilip}) can be proven when $A$ is any
topological space and $S(A)$ is any subset of $A^{\prime}$ such that the
weak*\ topology on $S(A)$ is metrizable. We will explore the relation between
\textrm{wu} and the algebraic structure on $A$ in later sections.
\end{remark}

Theorem (\ref{quasilip}) suggest the following definition:

\begin{definition}
A seminorm $L$ defined on a dense subset $\operatorname*{dom}L$ of the set of
self-adjoint elements of a separable C*-algebra $A$ and such that its dual
bounded-Lipschitz distance $d_{L}$ induces the restriction of the weak*
topology on the state space $S(A)$ of $A$ is called a quasi-Lip-norm.
\end{definition}

The terminology of quasi-Lip-norm is explained by the following simple result:

\begin{corollary}
Let $A$ be a unital separable C*-algebra and $L$ is a quasi-Lip-norm on $A$.
Then $L$ is a Lip-norm (as defined in \cite[Theorem 1.9]{Rieffel98a}) if and
only if the following two conditions hold:

\begin{enumerate}
\item $\{a\in\operatorname*{dom}L:L(a)=0\}=\mathbb{R}1$,

\item $\sup\left\{  d_{L}\left(  \varphi,\psi\right)  :\text{ }\varphi,\psi\in
S(A)\right\}  <\infty$.
\end{enumerate}
\end{corollary}

\begin{proof}
Apply Theorem \cite[Theorem 1.9]{Rieffel98a} with Theorem (\ref{quasilip}).
\end{proof}

On the other hand, whether $A$ is unital or not, the diameter of $S(A)$ for
$d_{L}$ is always 2. This may seem arbitrary (though see \cite[p. 2 and sec
1.7, pp. 24--27]{Weaver99}), however as far as the topology is concerned, we
could have chosen different bounds on $L$ and $\left\Vert .\right\Vert $ to
define a distance on $S(A)$ without much consequence for our purposes:

\begin{proposition}
\label{metricseq}Let $A$ be a separable C*-algebra and let $L$ a seminorm
defined on a dense subset $\operatorname*{dom}L$ of the self-adjoint part of
$A$. Let $\alpha,\beta>0$. The $(\alpha,\beta)$--bounded-Lipschitz metric
$d_{L,\alpha,\beta}$ defined on the state space $S(A)$ of $A$ by:
\[
d_{L,\alpha,\beta}(\varphi,\psi)=\sup\left\{  \left\vert \varphi
(a)-\psi(a)\right\vert :a\in\operatorname*{dom}L\text{, }L(a)\leq\beta\text{,
}\left\Vert a\right\Vert \leq\alpha\right\}
\]
is equivalent to the bounded-Lipschitz metric $d_{L}$. Consequently,
$d_{L,\alpha,\beta}$ metrizes the weak* topology on $S(A)$ if and only if $L $
is a quasi-Lip-norm.
\end{proposition}

\begin{proof}
Simply remark that if $\mathcal{B}_{L}^{\prime}=\left\{  a\in
\operatorname*{dom}L:L(a)\leq\beta\text{, }\left\Vert a\right\Vert \leq
\alpha\right\}  $ then $\left(  \max\{\alpha,\beta\}\right)  ^{-1}%
\mathcal{B}_{L}^{\prime}\subseteq\mathcal{B}_{L}\subseteq\left(  \min\left\{
\alpha,\beta\right\}  \right)  ^{-1}\mathcal{B}_{L}^{\prime}$.
\end{proof}

Informally, the family of metrics $(d_{L,\alpha,1})_{\alpha\in\mathbb{R}}$
describes the global metric geometry of $S(A)$ while any metric of the family
only describes the local metric properties of $S(A)$ (in particular one can
only recover $L$ from the whole family $(d_{L,\alpha,1})_{\alpha\in\mathbb{R}%
}$).

\bigskip To conclude this section, we propose a simple criterion on a seminorm
$L$ for the distance $d_{L}$ to be complete:

\begin{proposition}
\label{Gdelta} The state space $S(A)$ of a (separable) C*-algebra $A$ is a
$G_{\delta}$ in the weak* compact space $Q(A)$ of positive linear functionals
of norm at most $1$, so $S(A)$ is a Polish space: it is metrizable for the
weak*\ topology with a complete metric. Let $L$ be a quasi-Lip-norm on $A$,
and let $d_{L}$ be the bounded-Lipschitz distance induced by $L$ on $S(A)$.
Then $(S(A),d_{L})$ is a complete path-metric space if there exists an
approximate unit $(e_{n})_{n\in\mathbb{N}}$ in $A$ such that for all
$n\in\mathbb{N}$ we have $L(e_{n})\leq1$ and $\left\Vert e_{n}\right\Vert
\leq1$.
\end{proposition}

\begin{proof}
Simply observe $S(A)=%
{\displaystyle\bigcap\limits_{n=1}^{\infty}}
\left(  Q(A)\backslash\frac{n-1}{n}Q(A)\right)  $ where $\frac{n-1}{n}Q(A)$ is
trivially a weak* closed subset of $Q(A)$ to conclude that $S(A)$ is a
$G_{\delta}$. We have not used here that $A$ is separable.

Let $(e_{n})_{n\in\mathbb{N}}$ be an approximate unit of $A$ contained in
$\mathcal{B}_{L}$. Let $(\varphi_{n})_{n\in\mathbb{N}}$ be a Cauchy sequence
of states of $A$ for the bounded Lipschitz distance $d_{L}$. By linearity,
$(\varphi_{n}(a))_{n\in\mathbb{N}}$ is a Cauchy sequence for all
$a\in\operatorname*{dom}L$. Since $\operatorname*{dom}L$ is norm-dense in $A$,
and the sequence $(\varphi_{n})_{n\in\mathbb{N}}$ is norm-bounded by
assumption (it is a sequence of states), we conclude that for all $a\in A$ the
sequence $(\varphi_{n}(a))_{n\in\mathbb{N}}$ is a Cauchy sequence in
$\mathbb{C}$, so it converges to some number $\mu(a)$. It is immediate that
$\mu$, as the pointwise limit of continuous linear functionals on a Banach
space, is itself a continuous linear functional. It is also immediate that
$\mu$ is positive, and hence of norm $\lim_{n\rightarrow\infty}\mu(e_{n})$. On
the other hand, $(\varphi_{n})_{n\in\mathbb{N}}$ converges uniformly to $\mu$
on $\mathcal{B}_{L}$ by definition of $d_{L}$ (extended to $A^{\prime}$), so
we conclude $\lim_{n\rightarrow\infty}\mu(e_{n})=\lim_{k\rightarrow\infty}%
\lim_{n\rightarrow\infty}\varphi_{k}(e_{n})=1$. Hence $\mu\in S(A)$ and
$(\varphi_{n})_{n\in\mathbb{N}}$ converges to $\mu$ in $(S(A),d_{L})$.

Now it is easy to check that, if $\varphi,\psi\in S(A)$, then $\eta=\frac
{1}{2}\left(  \varphi+\psi\right)  \in S(A)$ as $S(A)$ is convex, and moreover
$d_{L}(\eta,\psi)=\frac{1}{2}d_{L}(\varphi,\psi)$. Since $S(A)$ is complete
for $d_{L}$, we deduce from \cite[Theorem 1.8 p. 7]{Gromov} that
$(S(A),d_{L})$ is a path metric space.
\end{proof}

Now, Theorem (\ref{quasilip}) relies upon a somewhat mysterious description of
\textrm{wu}. Our goal is now to reformulate Theorem\ (\ref{quasilip}) by using
a more explicit description of \textrm{wu} on bounded subsets of $A$. To
obtain this reformulation, we shall compare \textrm{wu} and the strict
topology on $A$.

\section{The Strict Topology and the weakly uniform topology}

\bigskip We let $A$ be a separable C*-algebra whose norm is denoted by
$\left\Vert .\right\Vert $ and whose state space $S(A)$ is endowed with the
weak* topology. In this section, we wish to compare our topology \textrm{wu}
to the strict topology. This is equivalent to investigate the joint continuity
of the multiplication of $A$ for \textrm{wu} on bounded sets, where a bounded
set will always mean bounded for the norm topology (we will justify this abuse
of language soon). Our approach consists in metrizing \textrm{wu} on bounded
subsets of $A$. Our main result is then that the strict topology and
\textrm{wu} agree on bounded subsets if and only the multiplication is jointly
continuous for \textrm{wu} on bounded subsets, which is not true in general
(but is true in the Abelian and unital case).

\subsection{First Elements of Comparison}

We start very simply with the trivial remark:

\begin{proposition}
\label{trivial}If $A$ is a unital C*-algebra, then the weakly-uniform topology
is the norm topology on $A$.
\end{proposition}

\begin{proof}
Since $A$ is unital, $S(A)$ $\in\mathfrak{S}$. Now $p_{S(A)}$ is equivalent to
$\left\Vert .\right\Vert $, and the seminorms $p_{K}$ for $K\in\mathfrak{S}$
are all dominated by $p_{S(A)}$, so \textrm{wu }is the norm topology.
\end{proof}

In the sequel, \textbf{we will not assume that }$A$ \textbf{is unital. }All of
the results which follow are valid in the unital case, but they are trivial
because of Proposition (\ref{trivial}). The reader may thus as well assume $A$
is not unital, and thus that $S(A)$ is not $\sigma(A^{\prime},A)$ compact.

The reader is probably already convinced that:

\begin{lemma}
\label{wusu copy(1)}\bigskip On any separable C*-algebra $A$, the
weakly-uniform topology \textrm{wu} is stronger than the weak topology
$\sigma(A,A^{\prime})$ on $A$, and in particular \textrm{wu} is Hausdorff.
\end{lemma}

\begin{proof}
Since any $\psi\in A^{\prime}$ is a linear combination of four states by
\cite[Secs. 3.1 and 3.2 pp. 41--46]{Pedersen79}, the smallest topology which
makes all states continuous is the weak* topology. Since $\left\{
\varphi\right\}  \in\mathfrak{S}$ for any $\varphi\in S(A)$, all the states of
$A$ are continuous for \textrm{wu}, hence our lemma.
\end{proof}

\bigskip In particular, a subset of $A$ is topologically bounded for any of
the \textrm{wu}, weak, strict or norm topologies if and only if it is bounded
for the norm topology by \cite[Theorem 3.18, p. 70]{Rudin91}. This justifies
that we shall simply call such a set bounded without reference to any topology.

We now offer a first comparison between \textrm{wu} and the strict topology.
This proposition differs from later comparisons by taking place on the whole
of $A$ rather than only on bounded subsets of $A$.

\begin{proposition}
\label{Cohen}The weakly-uniform topology is weaker than the strict topology on
any separable C*-algebra $A$.
\end{proposition}

\begin{proof}
For all $a,b\in A$ and $\mu\in A^{\prime}$ we define $\left(  a\cdot
\mu\right)  (b)=\mu(ab)$. The exterior composition law $\cdot$ gives
$A^{\prime}$ a Banach $A$--module structure. Let $K$ be a compact subset of
$S(A)$ for the weak*\ topology. Let $(e_{n})_{n\in\mathbb{N}}$ be any
approximate unit of $A$ which is also an increasing Abelian sequence of
positive elements of norm 1. Let us embed $A$ into its smallest unitalization
$\widetilde{A}$ of unit $1$. We then have, by Cauchy-Schwarz:
\begin{align*}
\left\Vert e_{n}\cdot\varphi-\varphi\right\Vert ^{2}  &  =\sup_{a\in
A,\left\Vert a\right\Vert \leq1}\left\vert e_{n}\cdot\varphi(a)-\varphi
(a)\right\vert ^{2}=\sup_{a\in A,\left\Vert a\right\Vert \leq1}\varphi
(a(e_{n}-1))^{2}\\
&  \leq\sup_{a\in A,\left\Vert a\right\Vert \leq1}\varphi(a^{\ast}%
a)\varphi\left(  \left(  e_{n}-1\right)  ^{2}\right)  \leq\varphi\left(
\left(  e_{n}-1\right)  ^{2}\right)  \text{.}%
\end{align*}

Since $\left(  e_{n}\right)  _{n\in\mathbb{N}}$ is Abelian, the sequence
$\left(  \left(  1-e_{n}\right)  ^{2}\right)  _{n\in\mathbb{N}}$ is
decreasing. As in Lemma (\ref{QuasiLipNeed}), for all $a\in A$ we define
$\widehat{a}:\varphi\in K\mapsto\varphi(a)$. Then $\left(  \widehat{\left(
1-e_{n}\right)  ^{2}}\right)  _{n\in\mathbb{N}}$ is a decreasing sequence of
continuous functions converging pointwise to $0$ in the compact $K$, hence
converges uniformly by Dini's theorem, so $\sup_{\varphi\in K}\left\Vert
e_{n}\cdot\varphi-\varphi\right\Vert \underset{n\rightarrow\infty
}{\longrightarrow}0$.

The set $K$ is bounded, so by \cite[Theorem 17.1, p. 114]{Doran79}, a
corollary the Cohen Factorization Theorem, there exists $c\in A$ and some
bounded subset $Y$ of $A^{\prime}$ such that $K=c\cdot Y$. Therefore, for all
$a\in A$:
\[
p_{K}(a)=\sup_{\varphi\in K}\left\vert \varphi(a)\right\vert =\sup_{\mu\in
Y}\left\vert \mu(ca)\right\vert \leq\left(  \sup_{\mu\in Y}\left\Vert
\mu\right\Vert \right)  \left\Vert ca\right\Vert \text{.}%
\]
Since $Y$ is bounded, the seminorm $p_{K}$ is strictly continuous.
\end{proof}

\subsection{Inductive Topologies}

We construct the bridge between the strict topology and the topology
\textrm{wu} in this section. This connection is established by introducing a
family of topologies which, at least on bounded subsets of $A$, will all agree
with \textrm{wu}, but are also very close in definition with the strict
topology. These topologies will be defined using sequences which, to ease
later estimates, we will allow to be in the enveloping Von Neumann algebra
$A^{\prime\prime}$ of $A$ (so that we can always choose sequences of
projections). Since the canonical map $A\rightarrow A^{\prime\prime}$ is a
*-monomorphism and hence an isometry, we still denote by $\left\Vert
.\right\Vert $ the norm on the enveloping Von Neumann algebra $A^{\prime
\prime}$ of $A$. Any state of $A$ is extended uniquely to a normal state of
$A^{\prime\prime}$. We define:

\begin{definition}
\label{ind copy(1)}Let $A$ be a separable C*-algebra whose bidual is denoted
by $A^{\prime\prime}$. Let $(e_{n})_{n\in\mathbb{N}}$ be a sequence in
$A^{\prime\prime}$. We assume moreover that $\left(  e_{n}\right)
_{n\in\mathbb{N}}$ satisfies the following condition:

\begin{condition}
\label{inductiveunits}%
\[
\forall\alpha\in\mathbb{N}\backslash\left\{  0\right\}  \ \ \ \forall
K\in\mathfrak{S\ \ \ \ }\lim_{n\rightarrow\infty}p_{K}\left(  \left(
1-e_{n}\right)  ^{\alpha}\right)  =0\text{.}%
\]

\end{condition}

We define on $A$ the seminorms $I_{e_{n}}:a\longmapsto\left\Vert e_{n}%
ae_{n}\right\Vert $ for all $n\in\mathbb{N}$. The locally convex topology on
$A$ generated by the seminorms $\left(  I_{e_{n}}\right)  _{n\in\mathbb{N}}$
is called the inductive topology along $(e_{n})_{n\in\mathbb{N}} $ and we will
denote it by \textrm{ind}$(e_{n})_{n\in\mathbb{N}}$.
\end{definition}

For Definition (\ref{ind copy(1)}) to be nonvacuous, we construct sequences
$\left(  e_{n}\right)  _{n\in\mathbb{N}}$ such that Condition
(\ref{inductiveunits}) holds. Since $A$ is separable, it contains a strictly
positive element, so we can define:

\begin{definition}
\label{approxunit}Let $A$ be a separable C*-algebra. Let $h\in A$ be a
strictly positive element such that $\left\Vert h\right\Vert =1$ and whose
spectrum in $A$ is denoted by $\sigma(h)$. Let $(f_{n})_{n\in\mathbb{N}}$ be a
sequence of nondecreasing functions on $\sigma(h)\cup\{0\}\subseteq
\lbrack0,1]$ such that for all $n\in\mathbb{N}$ we have $f_{n}(0)=0$ and
$f_{n}(1)=1$. We also assume that $(f_{n})_{n\in\mathbb{N}}$ converges
pointwise to $1$ on $\sigma(h)\backslash\{0\}$. The sequence $(f_{n}%
(h))_{n\in\mathbb{N}}$ of elements in the enveloping Von Neumann algebra
$A^{\prime\prime}$ of $A$ will be called a ($h$-)pseudo-spectral (ps-)
approximate unit of $A$.
\end{definition}

\begin{remark}
To distinguish ps-approximate units from the usual approximate units in $A$,
we shall call the latter continuous approximate units, and adopt the
convention that all continuous approximate units are positive and of norm
uniformly bounded by 1.
\end{remark}

Since any monotone function is a Borel function, ps-approximate units are
well-defined by the Borel functional calculus in $W^{\ast}(h)$ and are Abelian
sequences of positive elements of norm 1 in $A^{\prime\prime}$. We will prove
later in Proposition (\ref{approxunit0}) that ps-approximate units satisfy a
similar approximation property in \textrm{wu} to the defining approximation
property of continuous approximate units in the strict topology, yet for now
we do not need this fact. We now prove that there are many examples of
sequences verifying Condition (\ref{inductiveunits}):

\begin{lemma}
\label{Dini}Let $(e_{n})_{n\in\mathbb{N}}$ be an Abelian nondecreasing
continuous approximate unit of $A$ of positive elements of norm at most 1, or
a ps-approximate unit of $A$. Let $1$ denote the identity of $A^{\prime\prime
}$. Let $K\in\mathfrak{S}$. We have for all $\alpha\in\mathbb{N}%
\backslash\{0\}$:
\begin{equation}
\lim_{n\rightarrow\infty}\sup_{\varphi\in K}\left\vert \varphi\left(  \left(
1-e_{n}\right)  ^{\alpha}\right)  \right\vert =0\text{.}\label{Dini1}%
\end{equation}

\end{lemma}

\begin{proof}
We denote by $Q(A)$ the space of positive linear functionals of norm at most 1
on $A$, and we endow $Q(A)$ with the weak*\ topology. The space $Q(A)\ $is
often known as the quasi-state space of $A$ \cite{Pedersen79}. Let $a\in
A^{\prime\prime}$ be self-adjoint. The function $\widehat{a}$ defined on
$Q(A)$ by $\varphi\mapsto\varphi(a)$ is an affine Borel map vanishing at 0. By
\cite[Theorem 3.10.3, p. 70]{Pedersen79}, the map $\xi:a\in A^{\prime
\prime\operatorname*{sa}}\mapsto\widehat{a}$ is a linear order-preserving
isometry from $A^{\prime\prime\operatorname*{sa}}$ into the space
$B_{0}(Q(A))$ of real Borel affine functions on $Q(A)$ vanishing at 0, and the
image of $A^{\operatorname*{sa}}$ by $\xi$ is the subspace $A_{0}(Q(A))$ of
$B_{0}(Q(A))$ of real continuous affine functions on $Q(A)$ vanishing at 0.

We first recall the standard argument for Abelian continuous approximate
units. Let $K\in\mathfrak{S}$ and let $\alpha\in\mathbb{N}\backslash\{0\} $.
Since $(e_{n})_{n\in\mathbb{N}}$ is Abelian, the sequence $(\left(
1-e_{n}\right)  ^{\alpha})_{n\in\mathbb{N}}$ is also an Abelian nondecreasing
continuous approximate unit. The sequence $\left(  \widehat{\left(
1-e_{n}\right)  ^{\alpha}}\right)  _{n\in\mathbb{N}}$ is a nonincreasing
sequence of continuous affine functions on the compact $K$ whose limit is the
continuous function 0. Therefore, by Dini's theorem, $\left(  \widehat{\left(
1-e_{n}\right)  ^{\alpha}}\right)  _{n\in\mathbb{N}}$ converges uniformly on
$K$ to $0$, which proves (\ref{Dini1}) by definition. This part of the proof
did not require that $(e_{n})_{n\in\mathbb{N}}$ is Abelian when $\alpha=1$.

\bigskip We now turn to the case when $(e_{n})_{n\in\mathbb{N}}$ is a
ps-approximate unit. Let $K\in\mathfrak{S}$. Let $h$ be a strictly positive
element of norm 1 such that for all $n\in\mathbb{N}$, there exists a
nondecreasing function $f_{n}$ on $\sigma(h)\cup\{0\}$ such that $e_{n}%
=f_{n}(h)$ in $A^{\prime\prime}$. For any $\mu\in A^{\prime}$, we define
$\pi(\mu)$ to be its restriction to $C^{\ast}(h)$. By construction,
$\pi:A^{\prime}\rightarrow C^{\ast}(h)^{\prime}$ is a positive linear
surjection of norm 1. Moreover, since $C^{\ast}(h)$ contains a continuous
approximate unit since $h$ is strictly positive (e.g. $\left(  1-\left(
1-h\right)  ^{n}\right)  _{n\in\mathbb{N}}$), we have $\pi(\varphi)\in
S(C^{\ast}(h))$ for any state $\varphi\in S(A)$. Let $K^{\prime}%
=\pi(K)\subseteq S(C^{\ast}(h))$. As a norm-continuous linear function, $\pi$
is also a continuous map between $A^{\prime}$ and $C^{\ast}(h)^{\prime}$
endowed with their respective weak* topologies, so $K^{\prime}$ is weak*
compact in $S(C^{\ast}(h))$. Since for all $\alpha\in\mathbb{N}\backslash
\{0\}$ and all $n\in\mathbb{N}$, we have $\sup_{\varphi\in K}\left\vert
\varphi\left(  (1-e_{n})^{\alpha}\right)  \right\vert =\sup_{\varphi\in
K^{\prime}}\left\vert \varphi\left(  (1-e_{n})^{\alpha}\right)  \right\vert $,
we can now work solely in the Abelian C*-algebra $C^{\ast}(h)$. Since
$K^{\prime}$ is weak* compact, and $\sigma(h)\backslash\{0\}$ is a Polish
locally compact space, $K^{\prime}$ is uniformly tight. Let $\varepsilon>0$.
There exists a compact subset $k_{\varepsilon}$ of $\sigma(h)\backslash\{0\}$
such that $\sup\left\{  \varphi(\sigma(h)\backslash k_{\varepsilon}%
):\varphi\in K^{\prime}\right\}  \leq\varepsilon$, where we identify Radon
integrals with their Radon measures. Hence, for any $\alpha\in\mathbb{N}%
\backslash\{0\},n\in\mathbb{N}$ we have:%
\begin{equation}
\sup_{\varphi\in K^{\prime}}\left\vert \varphi\left(  \left(  1-f_{n}\right)
^{\alpha}\right)  \right\vert \leq\max\left\{  \varepsilon^{\alpha}%
,\sup_{\varphi\in K^{\prime}}\left(  \int_{k_{\varepsilon}}(1-f_{n})^{\alpha
}d\varphi\right)  \right\}  \text{.}\label{Dini0}%
\end{equation}

The sequence of nonincreasing functions $(1-f_{n})_{n\in\mathbb{N}}^{\alpha}$
converges pointwise on the compact subset $k_{\varepsilon}$ to the continuous
function 0 by assumption. Let $m=\min k_{\varepsilon}$ and $M=\max
k_{\varepsilon}$. Then we have, for all $n\in\mathbb{N}$ and all $x\in
k_{\varepsilon}$:%
\[
\left(  1-f_{n}\right)  ^{\alpha}(M)\leq\left(  1-f_{n}\right)  ^{\alpha
}(x)\leq\left(  1-f_{n}\right)  ^{\alpha}(m)
\]
so, since $\left(  1-f_{n}(m)\right)  _{n\in\mathbb{N}}^{\alpha}$ and $\left(
1-f_{n}(M)\right)  _{n\in\mathbb{N}}^{\alpha}$ converge to $0$ we deduce that
$(1-f_{n})_{n\in\mathbb{N}}^{\alpha}$ converges uniformly on $k_{\varepsilon}$
to 0 for any $\alpha\in\mathbb{N}\backslash\{0\}$. Therefore there exists
$N\in\mathbb{N}$ such that for all $n\geq N$ we have $\sup_{k_{\varepsilon}%
}(1-f_{n})^{\alpha}\leq\varepsilon^{\alpha}$. So, by (\ref{Dini0}), we have
$\sup_{\varphi\in K^{\prime}}\left\vert \varphi\left(  (1-e_{n})^{\alpha
}\right)  \right\vert \leq\varepsilon^{a}$.This concludes our proof since
$\varepsilon>0$ is arbitrary.
\end{proof}

\bigskip The main motivation for Definition (\ref{ind copy(1)}) is the
following proposition:

\begin{proposition}
\label{indwu}Let $A$ be a separable C*-algebra, and let $A^{\prime\prime}$ be
its bidual, which is the enveloping Von Neumann algebra of $A$. Let
$(e_{n})_{n\in\mathbb{N}}\in A^{\prime\prime}$ be a sequence satisfying
Condition (\ref{inductiveunits}). A bounded net converges for the
weakly-uniform topology \textrm{wu} if it converges in the inductive topology
\textrm{ind}$(e_{n})_{n\in\mathbb{N}}$. In other words, \textrm{ind}%
$(e_{n})_{n\in\mathbb{N}}$ is finer than \textrm{wu} on bounded subsets of $A
$.
\end{proposition}

\begin{proof}
\qquad Since any $a\in A$ can be written as a linear combination of the two
self-adjoint elements $\frac{1}{2}(a+a^{\ast})$ and $\frac{1}{2i}(a-a^{\ast}%
)$, and since both topologies make the addition, scalar multiplication and the
involution of $A$ continuous, we will restrict ourselves to nets of
self-adjoint elements of norm bounded by 1 and converging to 0.

Let now $(a_{\lambda})_{\lambda\in\Lambda}$ be a net in $A^{\operatorname*{sa}%
} $ converging to $0\in A$ for \textrm{ind}$(e_{n})_{n\in\mathbb{N}}$ and
uniformly bounded in norm by $1$. Let $K\in\mathfrak{S}$. Let $\varepsilon>0$.
By assumption, $\left(  e_{n}\right)  _{n\in\mathbb{N}}$ satisfies Condition
(\ref{inductiveunits}), so there exists $N\in\mathbb{N}$ such that, for all
$\varphi\in K$ and all $n\geq N$ we have $\varphi\left(  (1-e_{n})^{2}\right)
\leq\varepsilon^{2}$. Choose $n\geq N$ and $\varphi\in K$. Now:
\begin{align*}
\varphi(a_{\lambda}) &  =\varphi(e_{n}a_{\lambda}e_{n})+\varphi((1-e_{n}%
)a_{\lambda}e_{n})+\\
&  +\varphi(e_{n}a_{\lambda}(1-e_{n}))+\varphi((1-e_{n})a_{\lambda}%
(1-e_{n}))\text{.}%
\end{align*}
By the Cauchy-Schwarz inequality:
\begin{align*}
\left\vert \varphi((1-e_{n})a_{\lambda}e_{n})\right\vert  &  \leq
\sqrt[2]{\varphi((1-e_{n})^{2})}\sqrt[2]{\varphi(e_{n}a_{\lambda}^{\ast
}a_{\lambda}e_{n})}\\
&  \leq\sqrt[2]{\varphi((1-e_{n})^{2})}\leq\varepsilon\text{.}%
\end{align*}
Similarly $\left\vert \varphi(e_{n}a_{\lambda}(1-e_{n}))\right\vert
\leq\varepsilon$ and $\left\vert \varphi((1-e_{n})a_{\lambda}(1-e_{n}%
))\right\vert \leq\varepsilon$. Thus $\left\vert \varphi(a_{\lambda
})\right\vert \leq\left\vert \varphi(e_{n}a_{\lambda}e_{n})\right\vert
+3\varepsilon$, so:%
\[
p_{K}(a_{\lambda})\leq\sup_{\varphi\in K}\left\vert \varphi(e_{n}a_{\lambda
}e_{n})\right\vert +3\varepsilon\leq I_{e_{n}}(a_{\lambda})+3\varepsilon
\text{.}%
\]
Since the net $(a_{\lambda})_{\lambda\in\Lambda}$ converges to $0$ for
\textrm{ind}$(e_{n})_{n\in\mathbb{N}}$, there exists $\lambda_{0}\in\Lambda$
so that, for $\lambda\in\Lambda$ such that $\lambda\geq\lambda_{0}$ we have
$I_{e_{n}}(a_{\lambda})\leq\varepsilon$ and therefore $p_{K}(a_{\lambda}%
)\leq4\varepsilon$. So the net $(a_{\lambda})_{\lambda\in\Lambda}$ converges
to 0 for the \textrm{wu }topology.
\end{proof}

As a corollary of Proposition (\ref{indwu}), we have as claimed a first
element of comparison between the strict topology and \textrm{wu}:

\begin{corollary}
Let $A$ be a separable C*-algebra. The weakly-uniform topology \textrm{wu} on
any bounded subset $B$ of $A$ is weaker than the strict topology on $B$.
\end{corollary}

\begin{proof}
Let $(a_{\lambda})_{\lambda\in\Lambda}$ be a bounded net in $A$, converging
strictly to 0. Let $(e_{n})_{n\in\mathbb{N}}$ satisfy Condition
(\ref{inductiveunits}). Since the multiplication is (jointly) continuous on
$B$ for the strict topology, the bounded nets $(e_{n}a_{\lambda}%
e_{n})_{\lambda\in\Lambda}$ strictly converge to 0 for all $n\in\mathbb{N}$.
By definition of the strict topology we have $\lim_{\lambda\in\Lambda
}\left\Vert e_{n}a_{\lambda}e_{n}\right\Vert =0$ for all $n\in\mathbb{N}$.
Thus $(a_{\lambda})_{\lambda\in\Lambda}$ converges to $0$ in \textrm{wu} by
Proposition (\ref{indwu}).
\end{proof}

\bigskip Thanks to Proposition (\ref{indwu}), it is also possible to show that
the multiplication map is not jointly continuous on bounded sets for the
\textrm{wu} topology on the algebra $\mathcal{K}$ of compact operators on the
separable Hilbert space $\mathcal{H}$.

\begin{proposition}
\label{ctrexmp}In general, the multiplication is not jointly continuous on
bounded sets for the weakly-uniform topology \textrm{wu}. A counter-example is
given in the unit ball of the C*-algebra $\mathcal{K}$.
\end{proposition}

\begin{proof}
Let $(\zeta_{n})_{n\in\mathbb{N}}$ be an orthonormal basis for $\mathcal{H}$.
Denote the inner product of $\mathcal{H}$ by $\left\langle .,.\right\rangle $.
For each $n\in\mathbb{N}$ we define the compact operator $S_{n}:\xi
\in\mathcal{H}\longmapsto\left\langle \xi,\zeta_{0}\right\rangle \zeta_{n}$.
It is easy to see that, for all $n\in\mathbb{N}$, we have $S_{n}^{\ast}%
(\xi)=\left\langle \xi,\zeta_{n}\right\rangle \zeta_{0}$ and thus $S_{n}%
^{\ast}S_{n}(\xi)=\left\langle \xi,\zeta_{0}\right\rangle \zeta_{0}$. Of
course, $(S_{n})_{n\in\mathbb{N}}$, $(S_{n}^{\ast})_{n\in\mathbb{N}}$ and
$(S_{n}^{\ast}S_{n})_{n\in\mathbb{N}}$ are bounded sequences. For all
$n\in\mathbb{N}$, we let $P_{n}$ be the projection in $\mathcal{H}$ on the
span of $\{\zeta_{0},\ldots,\zeta_{n}\}$. The sequence $(P_{n})_{n\in
\mathbb{N}}$ is an (continuous) increasing positive Abelian approximate unit
for $\mathcal{K}$ and of course $\left\Vert P_{n}\right\Vert =1$ for all
$n\in\mathbb{N}$. In addition, we can check easily that, for all
$n\in\mathbb{N}$ and for all $k>n$ we have $P_{n}S_{k}P_{n}=0$ and $P_{n}%
S_{k}^{\ast}P_{n}=0$. Hence, $(S_{n})_{n\in\mathbb{N}}$ and $(S_{n}^{\ast
})_{n\in\mathbb{N}}$ converge to 0 in \textrm{ind}$(P_{n})_{n\in\mathbb{N}}$,
and therefore in the \textrm{wu} topology by Proposition (\ref{indwu}). Yet
$(S_{n}^{\ast}S_{n})_{n\in\mathbb{N}}$ is constant and nonzero.
\end{proof}

\subsection{Algebraic Description of the Weakly-Uniform Topology}

We now use a specific type of ps-approximate units to give our first algebraic
description of the topology \textrm{wu}:

\begin{definition}
\label{hspectral}Let $A$ be a separable C*-algebra. Let $h\in A$ be a strictly
positive element such that $\left\Vert h\right\Vert =1$. A sequence $\left(
h_{n}\right)  _{n\in\mathbb{N}}$ is a $h$-spectral approximate unit when there
exists a decreasing sequence $\left(  \alpha_{n}\right)  _{n\in\mathbb{N}}$ in
$]0,1]$ converging to 0 such that for each $n\in\mathbb{N}$ we have
$h_{n}=\chi_{\lbrack\alpha_{n},1]}(h)$, where $\chi_{E}$ is the indicator
function of any subset $E$ in $[0,1]$. In particular, by definition, $\left(
h_{n}\right)  _{n\in\mathbb{N}}$ is a ps-approximate unit of projections in
$A^{\prime\prime}$.
\end{definition}

\begin{theorem}
\label{wu_spectral}Let $A$ be a separable C*-algebra, and let $(a_{\lambda
})_{\lambda\in\Lambda}$ be a bounded net in $A$. Let $h\in A$ be a strictly
positive element such that $\left\Vert h\right\Vert =1$, and $(h_{n}%
)_{n\in\mathbb{N}}$ be a $h$-spectral approximate unit for $A$. Then
$(a_{\lambda})_{\lambda\in\Lambda}$ converges to 0 in the weakly-uniform
topology \textrm{wu} if, and only if it converges to 0 in the inductive
topology \textrm{ind}$\left(  h_{n}\right)  _{n\in\mathbb{N}}$.
\end{theorem}

\begin{proof}
Since $h$-spectral approximate units are ps-approximate unit, by Lemma
(\ref{Dini}) and Proposition (\ref{indwu}), the topology \textrm{ind}$\left(
h_{n}\right)  _{n\in\mathbb{N}}$ is stronger than \textrm{wu} on any bounded
set. We now turn to the converse inclusion.

Let $n\in\mathbb{N}$. Note that $h_{n}\in A^{\prime\prime}$ is a projection.
For all $t\in\sigma(h)$ we set $g_{n}(t)=1$ if $t\geq\alpha_{n}$ and
$g_{n}(t)=\alpha_{n}^{-1}t$ otherwise, and we set $e_{n}=g_{n}(h)\in C^{\ast
}(h)$. Obviously, $e_{n}h_{n}=h_{n}e_{n}=h_{n}$. We check easily that, since
$h_{n}$ is a projection, if $a\in A$ then we have $\left(  h_{n}ah_{n}\right)
\left(  h_{n}e_{n}h_{n}\right)  =h_{n}ah_{n}e_{n}h_{n}=h_{n}ah_{n}$ and
similarly $\left(  h_{n}e_{n}h_{n}\right)  \left(  h_{n}ah_{n}\right)
=h_{n}ah_{n}$, so that $h_{n}\in h_{n}Ah_{n}$ and $h_{n}$ is the unit of
$h_{n}Ah_{n}$ (even though $A$ has no unit in general). Let $\varphi\in
S(h_{n}Ah_{n})$. Define $\iota_{n}(\varphi)(a)=\varphi(h_{n}ah_{n})$ for all
$a\in A$. The map $\iota_{n}(\varphi)$ is obviously a positive linear
functional, since if $a$ is positive then $h_{n}ah_{n}$ is positive as well,
and $\left\Vert \iota(\varphi)\right\Vert \leq1$. On the other, $\iota
_{n}(\varphi)(e_{n})=1$, and since $\left\Vert e_{n}\right\Vert =1$ we deduce
that $\iota_{n}(\varphi)$ is indeed a state of $A$. The map $\iota_{n}%
:S(h_{n}Ah_{n})\rightarrow S(A)$ is thus an affine map which is obviously
continuous when the state spaces are endowed with the weak* topology (it is
also immediate that $\iota$ is injective). Let $K_{n}$ be the range of
$\iota_{n}$: by continuity of $\iota_{n}$, the set $K_{n}$ is weak*-compact in
$S(A)$.

Let now $\left(  a_{\lambda}\right)  _{\lambda\in\Lambda}$ be a bounded net in
$A$ converging to $0$ for \textrm{wu}. For $n\in\mathbb{N}$ and $\lambda
\in\Lambda$ we have:%
\begin{align*}
\left\Vert h_{n}a_{\lambda}h_{n}\right\Vert  &  =\sup_{\varphi\in
S(h_{n}Ah_{n})}\varphi(h_{n}a_{\lambda}h_{n})=\sup_{\varphi\in S(h_{n}Ah_{n}%
)}\iota_{n}(\varphi)(a_{\lambda})\\
&  =\sup_{\varphi\in K_{n}}\varphi(a_{\lambda})=p_{K_{n}}(a_{\lambda})\text{.}%
\end{align*}
Hence, since $a_{\lambda}\underset{\lambda\in\Lambda}{\overset{\mathrm{wu}%
}{\longrightarrow}}0$ we have $\lim_{\lambda\in\Lambda}\left\Vert
h_{n}a_{\lambda}h_{n}\right\Vert =\lim_{\lambda\in\Lambda}p_{K_{n}}%
(a_{\lambda})=0$ for all $n\in\mathbb{N}$ so $\left(  a_{\lambda}\right)
_{\lambda\in\Lambda}$ converges to $0$ in the topology \textrm{ind}$\left(
h_{n}\right)  _{n\in\mathbb{N}}$. This concludes the proof of our Theorem.
\end{proof}

\bigskip Theorem (\ref{wu_spectral}) gives us our first insight into a simpler
description for the \textrm{wu}--topology. It is still a bit difficult to
handle a spectral approximate unit, so we shall simplify this result one step
further. The following lemma extends mildly on \cite[Lemma 3.2.6]%
{Wegge-Olsen92}:

\begin{lemma}
\label{strict}Let $A$ be a separable C*-algebra, $h\in A$ a strictly positive
element. Let $(a_{\lambda})_{\lambda\in\Lambda}$ be a bounded net in $A$. Then%
\begin{align}
\lim_{\lambda\in\Lambda}\left\Vert ha_{\lambda}\right\Vert  &  =0\text{ if and
only if }\lim_{\lambda\in\Lambda}\left\Vert ba_{\lambda}\right\Vert =0\text{
for all }b\in A\text{,}\label{strict0}\\
\lim_{\lambda\in\Lambda}\left\Vert a_{\lambda}h\right\Vert  &  =0\text{ if and
only if }\lim_{\lambda\in\Lambda}\left\Vert a_{\lambda}c\right\Vert =0\text{
for all }c\in A\text{,}\label{strict1}\\
\lim_{\lambda\in\Lambda}\left\Vert ha_{\lambda}h\right\Vert  &  =0\text{ if
and only if }\lim_{\lambda\in\Lambda}\left\Vert ba_{\lambda}c\right\Vert
=0\text{ for all }b,c\in A\text{.}\label{strict2}%
\end{align}

In particular, the strict topology on any bounded subset $h$ of $A$ agree with
the topology defined on $h$ by the two seminorms
\[
a\longmapsto\left\Vert ha\right\Vert \text{ and }a\longmapsto\left\Vert
ah\right\Vert \text{.}%
\]

\end{lemma}

\begin{proof}
The sufficient conditions above are all trivial, by setting $b=c=h,$ and the
necessary conditions of (\ref{strict0}) and (\ref{strict1}) are found in
\cite[Lemma 3.2.6]{Wegge-Olsen92}, together with the metrizability of the
strict topology on bounded subsets of $A$. Now, let $\left(  b,c\right)  \in
A$ and assume $\lim_{\lambda\in\Lambda}\left\Vert ha_{\lambda}h\right\Vert
=0$. In particular, the net $\left(  a_{\lambda}h\right)  _{\lambda\in\Lambda
}$ converges to $0$ for the seminorm $a\mapsto\left\Vert ha\right\Vert $ and
thus, by (\ref{strict0}), we conclude that for all $b\in A$ we have
$\lim_{\lambda\in\Lambda}\left\Vert ba_{\lambda}h\right\Vert =0$. By
(\ref{strict1}), since the net $\left(  ba_{\lambda}\right)  _{\lambda
\in\Lambda}$ converges to $0$ for the seminorm $a\mapsto\left\Vert
ah\right\Vert $, we conclude that $\lim_{\lambda\in\Lambda}\left\Vert
ba_{\lambda}c\right\Vert =0$ as required.
\end{proof}

To describe a metric for \textrm{wu} on bounded sets, we need to further our
understanding of ps-approximate units. Ps-approximate units are meant to have
a similar approximation property in \textrm{wu} as the continuous approximate
units have for the strict topology. We shall now establish this, using the
following easy lemma as a tool:

\begin{lemma}
\label{AbelianCompact}Let $X$ be a locally compact Polish space. Let $k$ be a
compact subset of $X$. Then the set $K=\{\mu\in S(C_{0}(X)):\mu(k)=1\}$ of
probability measures supported on $k$ is a weak* compact subset of
$S(C_{0}(X))$ where $C_{0}(X)$ the is C*-algebra of continuous functions on
the one-point compactification on $X$ vanishing at infinity, and we identify
the continuous linear functionals on $C_{0}(X)$ with the bounded Radon
measures on $X$ by the Radon-Riesz theorem.
\end{lemma}

\begin{proof}
Let $(\mu_{n})_{n\in\mathbb{N}}$ be a sequence in $K$. The set $K$ is
uniformly tight by construction by \cite[Theorem 11.5.4., p.404]{Dudley}, so
$(\mu_{n})_{n\in\mathbb{N}}$ admits a subsequence $(\mu_{m(n)})_{n\in
\mathbb{N}}$ which converges to $\mu\in S(C_{0}(X))$ for the weak* topology.
Now, assume that $\mu(k)<1$. Then there exists an open set $U$ in $X\backslash
k$ so that $\mu(U)>0$. Hence there exists a nonzero positive compactly
supported continuous function $f$ on $U$ and some $\varepsilon>0$ such that
$\mu(f)=\varepsilon>0$. We extend $f$ to $X$ by setting $f(X\backslash U)=0$.
Now, $f\in C_{0}(X)$ and therefore $\lim_{n\rightarrow\infty}\mu_{m(n)}%
(f)=\mu(f)=\varepsilon$. Yet, since $f(k)=\{0\}$ by construction, we have
$\mu_{n}(f)=0$ for all $n\in\mathbb{N}$, so we reached a contradiction. Hence
$\mu(k)=1$ and so $\mu\in K$ which concludes the proof of our lemma, since $X$
is separable so the weak*$\ $topology of $S(C_{0}(X))\ $is metrizable.
\end{proof}

We now prove the approximation property of ps-approximate units:

\begin{lemma}
\label{approxunit0}Let $A$ be a separable C*-algebra. Let $h\in A$ be a
strictly positive element and let $(e_{n})_{n\in\mathbb{N}}$ be an
$h$-ps-approximate unit. Then for all $a\in A$, the sequences $(ae_{n}%
)_{n\in\mathbb{N}}$ and $(e_{n}a)_{n\in\mathbb{N}}$ converge to $a$ for the
weakly uniform topology \textrm{wu}. Moreover, for all $a\in C^{\ast}(h)$:%
\[
\lim_{n\rightarrow\infty}\left\Vert ae_{n}-a\right\Vert =\lim_{n\rightarrow
\infty}\left\Vert e_{n}a-a\right\Vert =0\text{.}%
\]

\end{lemma}

\begin{proof}
Let $h\in A$ be a strictly positive element such that $\left\Vert h\right\Vert
=1$ and $f_{n}(h)=e_{n}$ where $f_{n}$ is a nondecreasing function on
$\sigma(h)\cup\{0\}$ such that $f_{n}(0)=0$, and $\left\Vert f_{n}\right\Vert
\leq1$ for all $n\in\mathbb{N}$. Let $a\in A$. Let $K\in\mathfrak{S}$. For any
$\varphi\in K$ we have, by Cauchy-Schwarz:%
\begin{equation}
\left\vert \varphi(ae_{n}-a)\right\vert ^{2}\leq\varphi(a^{\ast}%
a)\varphi\left(  \left(  1-e_{n}\right)  ^{2}\right) \label{approxunit0.2}%
\end{equation}
where, as usual, 1 is the unit of $A^{\prime\prime}$. Hence, taking the
supremum over $\varphi\in K$ in (\ref{approxunit0.2}) we get $p_{K}%
(ae_{n}-a)^{2}\leq\left\Vert a\right\Vert ^{2}p_{K}\left(  \left(
1-e_{n}\right)  ^{2}\right)  $.

By Lemma (\ref{Dini}), we have $\lim_{n\rightarrow\infty}p_{K}\left(  \left(
1-e_{n}\right)  ^{2}\right)  =0$, so we conclude that $p_{K}\left(
ae_{n}-a\right)  \underset{n\rightarrow\infty}{\longrightarrow}0$ for all
$K\in\mathfrak{S}$. The same reasoning holds for the sequence $(e_{n}%
a)_{n\in\mathbb{N}}$.

Now, let $g$ be a continuous function over $\sigma(h)\cup\{0\}$ such that
$g(0)=0$. Since $g(0)=0$ and $g$ is continuous on $\sigma(h)\subseteq
\lbrack0,1]$, there exists $z_{0}\in\sigma(h)\backslash\{0\}$ such that for
all $z\in\sigma(h)$ such that $z<z_{0}$ we have $\left\vert g(z)\right\vert
\leq\frac{1}{2}\varepsilon$. Let $k=\left\{  z\in\sigma(h):z\geq
z_{0}\right\}  $. By construction, $k$ is a closed subset of the compact
$\sigma(h)$, and we note also $k\subseteq\sigma(h)\backslash\{0\} $ . By Lemma
(\ref{AbelianCompact}), the set $K=\left\{  \mu\in S(C_{0}(\sigma
(h)):\mu(k)=1\right\}  $ is a weak* compact subset of the state space of
$C_{0}(\sigma(h)\backslash\{0\})$. Now, denoting by $\left\Vert .\right\Vert
_{\infty}$ the supremum norm on $L^{\infty}(\sigma(h)\cup\{0\})$, we have:%
\begin{equation}
\left\Vert f_{n}g-g\right\Vert _{\infty}\leq\sup_{x\in k}\left\{  \left\vert
f_{n}(x)g(x)-g(x)\right\vert ,\varepsilon\right\} \label{approxunit0.0}%
\end{equation}
since $\left\vert f_{n}(x)\right\vert \leq1$ by assumption, and since for
$x\not \in k$ we have $\left\vert g(x)\right\vert \leq\frac{1}{2}\varepsilon$.
Now, we have%
\begin{equation}
\sup_{x\in k}\left\vert f_{n}(x)g(x)-g(x)\right\vert \leq\left\Vert
g\right\Vert _{\infty}\sup_{x\in k}\left\vert f_{n}(x)-1\right\vert
\text{.}\label{approxunit0.1}%
\end{equation}
The sequence $\left(  1-f_{n}\right)  _{n\in\mathbb{N}}$ of nonincreasing
functions converges on the compact $k$ pointwise, hence uniformly to $0$, as
in Lemma (\ref{Dini}). Therefore, there exists $N\in\mathbb{N}$ such that
$\left\Vert g\right\Vert _{\infty}\sup_{x\in k}\left\vert f_{n}%
(x)-1\right\vert \leq\varepsilon$ for all $n\geq N$. Thus by
(\ref{approxunit0.0}) and (\ref{approxunit0.1}) we have for all $n\geq N$ that
$\left\Vert f_{n}g-g\right\Vert _{\infty}\leq\varepsilon$, so $\lim
_{n\rightarrow\infty}\left\Vert f_{n}g-g\right\Vert _{\infty}=0$ for all $g\in
C_{0}(\sigma(h)\backslash\{0\})$. Consequently for all $a\in C^{\ast}(h)$ we
have $\lim_{n\rightarrow\infty}\left\Vert f_{n}(h)a-a\right\Vert =0$.
\end{proof}

\bigskip We can now obtain the main theorem about \textrm{wu}, which provides
a simple description of this topology on bounded subsets of $A$ as a metric
topology:{}

\begin{theorem}
\label{wunorm}Let $A$ be a separable C*-algebra. Let $h\in A$ be a strictly
positive element. Let $(a_{\lambda})_{\lambda\in\Lambda}$ be a bounded net in
$A$. Then $(a_{\lambda})_{\lambda\in\Lambda}$ converges to 0 in the
weakly-uniform topology \textrm{wu} if, and only if $\lim_{\lambda\in\Lambda
}\left\Vert ha_{\lambda}h\right\Vert =0$. Equivalently, on any bounded set
$B\subseteq A$, the distance $\Delta:a,b\in B\mapsto\left\Vert
h(b-a)h\right\Vert $ metrizes \textrm{wu}.
\end{theorem}

\begin{proof}
Let us first assume that the bounded net $(a_{\lambda})_{\lambda\in\Lambda}$
converging to 0 in the topology \textrm{wu}. To simplify, assume $\left\Vert
a_{\lambda}\right\Vert \leq1$ for all $\lambda\in\Lambda$. Now, let $h\in A$
strictly positive, and assume without loss of generality that $\left\Vert
h\right\Vert =1$. Let $\left(  h_{n}\right)  _{n\in\mathbb{N}}$ be a
$h$--spectral approximate unit. Then, for all $n\in\mathbb{N}$, by Theorem
(\ref{wu_spectral}) we have $\lim_{\lambda\in\Lambda}\left\Vert h_{n}%
a_{\lambda}h_{n}\right\Vert =0$. Let $\varepsilon>0$. There exists
$N\in\mathbb{N}$ such that, for all $n\geq N$, we have $\left\Vert
h-hh_{n}\right\Vert \leq\frac{\varepsilon}{2}$ and $\left\Vert h-h_{n}%
h\right\Vert \leq\frac{\varepsilon}{2}$ by Lemma (\ref{approxunit0}). We then
have:
\begin{align*}
\left\Vert ha_{\lambda}h\right\Vert  &  \leq\left\Vert hh_{n}a_{\lambda
}h\right\Vert +\left\Vert ha_{\lambda}h-hh_{n}a_{\lambda}h\right\Vert \\
&  \leq\left\Vert h\right\Vert \left\Vert h_{n}a_{\lambda}h\right\Vert
+\left\Vert h-hh_{n}\right\Vert \left\Vert a_{\lambda}h\right\Vert
\leq\left\Vert h_{n}a_{\lambda}h\right\Vert +\frac{\varepsilon}{2}\\
&  \leq\left\Vert h_{n}a_{\lambda}h_{n}h\right\Vert +\left\Vert h_{n}%
a_{\lambda}h-h_{n}a_{\lambda}h_{n}h\right\Vert +\frac{\varepsilon}{2}\\
&  \leq\left\Vert h_{n}a_{\lambda}h_{n}\right\Vert \left\Vert h\right\Vert
+\left\Vert h_{n}a_{\lambda}\right\Vert \left\Vert h-h_{n}h\right\Vert
+\frac{\varepsilon}{2}\leq\left\Vert h_{n}a_{\lambda}h_{n}\right\Vert
+\varepsilon\text{.}%
\end{align*}
Hence $0\leq\overline{\lim_{\lambda\in\Lambda}}\left\Vert ha_{\lambda
}h\right\Vert \leq\varepsilon$, and since $\varepsilon>0$ is arbitrary, we
conclude that $\overline{\lim}_{\lambda\in\Lambda}\left\Vert ha_{\lambda
}h\right\Vert =0$ and subsequently $\lim_{\lambda\in\Lambda}\left\Vert
ha_{\lambda}h\right\Vert =0$.

Conversely, we now assume that $\lim_{\lambda\in\Lambda}\left\Vert
ha_{\lambda}h\right\Vert =0$, with the same notations as above. By the Lemma
(\ref{strict}), for any elements $b,c\in A$, we have that $\lim_{\lambda
\in\Lambda}\left\Vert ba_{\lambda}c\right\Vert =0$. This implies that
$(a_{\lambda})_{\lambda\in\Lambda}$ converges to 0 in \textrm{ind}%
$(e_{n})_{n\in\mathbb{N}}$ for any continuous, Abelian approximate unit
$(e_{n})_{n\in\mathbb{N}}$ in $A$. This implies \textrm{wu} convergence, as
the nets are bounded, by Proposition (\ref{indwu}). This concludes our theorem.
\end{proof}

\begin{corollary}
\label{wunorm0}Let $A$ be a separable C*-algebra. The weakly uniform topology
restricted to any bounded subsets $A$ is generated by the set of seminorms
$\left(  a\in A\longmapsto\left\Vert bac\right\Vert \right)  _{b,c\in A}$.
\end{corollary}

In the next section, we try to return to the state space picture and
understand what differs between the strict topology and \textrm{wu}.

\subsection{The Strong Uniform Topology}

We now introduce the topology \textrm{su}, which is modelled after
\textrm{wu}, but constructed with the joint continuity of the multiplication
on bounded subsets of $A$ in mind as well. The question of the \textrm{wu}%
-joint continuity of the multiplication on bounded subsets of $A$ is
essentially the question of when \textrm{wu }agrees with \textrm{su }on
bounded subsets, as we shall see in Theorem (\ref{sustrict}) and
(\ref{wustrict}). It turns out that \textrm{su} restricted to bounded subsets
of $A$ is the strict topology.

We now introduce our topology \textrm{su }on $A$:

\begin{definition}
Let $A$ be a separable C*-algebra, and let $\mathfrak{S}$ be the set of all
weak* compact subsets of the state space $S(A)$ of $A$. For any $K\in
\mathfrak{S}$, we define the seminorms $q_{K}$ by:
\[
q_{K}(a)=\sup_{\varphi\in K}\left(  \sqrt[2]{\varphi(a^{\ast}a)}%
,\sqrt[2]{\varphi(aa^{\ast})}\right)
\]
for all $a\in A$. The locally convex topology on $A$ generated by the family
$(q_{K})_{K\in\mathfrak{S}}$ of seminorms will be called the strongly uniform
topology, or $\mathrm{su}$ topology.
\end{definition}

Using Theorem (\ref{wunorm}) it is immediate that:

\begin{theorem}
\label{sustrict}Let $A$ be a separable C*-algebra. The strict topology and the
strong-uniform topology \textrm{su} agree on the bounded subsets of $A$.
\end{theorem}

\begin{proof}
Let $(a_{\lambda})_{\lambda\in\Lambda}$ be a bounded net is $A$. Now
$(a_{\lambda})_{\lambda\in\Lambda}\overset{\mathrm{su}}{\underset{\lambda
\in\Lambda}{\longrightarrow}}0$ if, and only if $(a_{\lambda}^{\ast}%
a_{\lambda})_{\lambda\in\Lambda}\overset{\mathrm{wu}}{\underset{\lambda
\in\Lambda}{\longrightarrow}}0$ and $(a_{\lambda}a_{\lambda}^{\ast}%
)_{\lambda\in\Lambda}\overset{\mathrm{wu}}{\underset{\lambda\in\Lambda
}{\longrightarrow}}0$ by definition of \textrm{su}. This in turn is
equivalent, by Theorem (\ref{wunorm}), to:%
\[
0=\lim_{\lambda\in\Lambda}\left\Vert ha_{\lambda}^{\ast}a_{\lambda
}h\right\Vert =\lim_{\lambda\in\Lambda}\left\Vert a_{\lambda}h\right\Vert
^{2}\text{ and }0=\lim_{\lambda\in\Lambda}\left\Vert ha_{\lambda}a_{\lambda
}^{\ast}h\right\Vert =\lim_{\lambda\in\Lambda}\left\Vert ha_{\lambda
}\right\Vert ^{2}%
\]
which, by Lemma (\ref{strict}), is equivalent to the strict convergence of
$(a_{\lambda})_{\lambda\in\Lambda}$ to 0.
\end{proof}

We can now prove the following necessary and sufficient condition on
\textrm{wu} for the multiplication to be \textrm{wu}-jointly continuous:

\begin{theorem}
\label{wustrict}Let $A$ be a separable C*-algebra. The multiplication of $A$
is jointly continuous for the weakly-uniform topology \textrm{wu }on any
bounded subset $B$ of $A$ if and only if \textrm{wu }and the strict topology
agree on $B$.
\end{theorem}

\begin{proof}
Let $B$ be a bounded subset of $A$. Let $(a_{\lambda})_{\lambda\in\Lambda}$ be
a net in $B$ converging for the topology \textrm{wu} to $a\in B$. Assume that
the multiplication is jointly \textrm{wu}-continuous on $B$. Then the nets
$(a^{\ast}a_{\lambda})_{\lambda\in\Lambda}$, $(a_{\lambda}^{\ast}%
a)_{\lambda\in\Lambda}$ and $(a_{\lambda}^{\ast}a_{\lambda})_{\lambda
\in\Lambda} $ all converge to $a^{\ast}a$ for \textrm{wu}, while the nets
$(aa_{\lambda}^{\ast})_{\lambda\in\Lambda}$, $(a_{\lambda}a^{\ast}%
)_{\lambda\in\Lambda}$ and $(a_{\lambda}a_{\lambda}^{\ast})_{\lambda\in
\Lambda}$ all converge to $aa^{\ast}$ for \textrm{wu. }We deduce that for any
$K\in\mathfrak{S}$, we have $\lim_{\lambda\in\Lambda}\sup_{\varphi\in
K}\varphi\left(  \left(  a_{\lambda}-a\right)  ^{\ast}\left(  a_{\lambda
}-a\right)  \right)  =0$ and $\lim_{\lambda\in\Lambda}\sup_{\varphi\in
K}\varphi\left(  \left(  a_{\lambda}-a\right)  \left(  a_{\lambda}-a\right)
^{\ast}\right)  =0$. Hence $(a_{\lambda})_{\lambda\in\Lambda}$ converges to $a
$ for \textrm{su}. By Theorem (\ref{sustrict}), \textrm{su} is the strict
topology when restricted to bounded sets, so $(a_{\lambda})_{\lambda\in
\Lambda}$ converges strictly to $a$: hence \textrm{wu} is finer than the
strict topology on bounded sets. Since by Proposition (\ref{Cohen}), the
strict topology is finer than \textrm{wu}, this conclude the necessary
condition. The sufficient condition is trivial.
\end{proof}

We already encountered in Proposition (\ref{ctrexmp}) a situation when
\textrm{wu} is strictly weaker than the strict topology on a bounded set. We
also know that trivially, \textrm{wu} and the strict topology agree (with the
norm topology) when $A$ is unital. Another situation when \textrm{wu} and the
strict topology agree on bounded sets is when $A$ contains a central, strictly
positive element:

\begin{proposition}
\label{strictwuAbelian}Let $A$ be a separable C*-algebra with a central
strictly positive element. Then the weakly-uniform topology \textrm{wu} and
the strict topology agree on bounded subsets of $A$.
\end{proposition}

\begin{proof}
Let $h\in A$ be a central strictly positive. On any bounded subset $B$ of $A$,
the topology \textrm{wu} is generated by the seminorm $a\mapsto\left\Vert
hah\right\Vert =\left\Vert h^{2}a\right\Vert =\left\Vert ah^{2}\right\Vert $,
and since $h^{2}$ is also strictly positive, this seminorm also generates the
strict topology on $B$.
\end{proof}

\begin{corollary}
If $A$ is an Abelian separable C*-algebra then the strict topology and the
\textrm{wu} topology agree on bounded subsets of $A$.
\end{corollary}

\begin{example}
For a non-Abelian example, consider $A=C_{0}(\mathbb{R})\otimes M_{2}$, where
$M_{2}$ is the algebra of $2\times2$ complex matrices. It is neither unital
nor Abelian. On the other hand, if $f:x\in\mathbb{R}\longmapsto\exp(-x^{2})$
then it is straightforward to check that $f\otimes I_{2}$, where $I_{2}$ is
the identity of $M_{2}$, is a strictly positive central element in $A$, hence
Proposition (\ref{strictwuAbelian}) applies to $A$.
\end{example}

Proposition (\ref{strictwuAbelian}) illustrates that, on bounded subsets of
$A$, the distinction between \textrm{wu} and the strict topology only appears
when, informally, $A$ has a high degree of noncommutativity, and that the
nonunital Abelian case may be sometimes misleading in the further study of
locally compact quantum metric spaces.

\section{Application to Quasi-Lip-norms, Lip-norms and finite diameter quantum
metric spaces}

We can now rephrase Theorem\ (\ref{quasilip}) using Theorem (\ref{wunorm}) and
Corollary (\ref{wunorm0}) to provide a satisfactory answer to Problem
(\ref{Pb}):

\begin{theorem}
\label{quasiLipalg}Let $A$ be a separable C*-algebra. Let $L$ be a seminorm on
a dense subspace $\operatorname*{dom}L$ of the self-adjoint part
$A^{\operatorname*{sa}}$ of a C*-algebra $A$. Let $\mathcal{B}_{L}=\left\{
a\in\operatorname*{dom}L:L(a)\leq1,\left\Vert a\right\Vert \leq1\right\}  $.
Then the following are equivalent:

\begin{itemize}
\item The seminorm $L$ is a quasi-Lip-norm (namely $L$ solves Problem
(\ref{Pb})),

\item For all $b,c\in A$, the set $b\mathcal{B}_{L}c$ is norm precompact in
$A$,

\item There exists a sequence $\left(  e_{n}\right)  _{n\in\mathbb{N}}$
satisfying Condition (\ref{inductiveunits}) such that $e_{n}\mathcal{B}%
_{L}e_{n}$ is norm totally bounded for all $n\in\mathbb{N}$.

\item There exists a strictly positive element $h\in A$ such that
$h\mathcal{B}_{L}h$ is norm-precompact in $A$.
\end{itemize}
\end{theorem}

\bigskip Our setting is particularly satisfactory with finite diameter locally
compact quantum metric spaces, since we then get:

\begin{proposition}
Let $A$ be a separable C*-algebra. Let $L$ be a quasi-Lip-norm on $A$ and
denote by $\kappa_{L}$ the distance induced on the state space $S(A)$ of $A$
by the dual of the seminorm $L$. Suppose $\kappa_{L}$ gives $S(A)$ a finite
diameter. Then $\kappa_{L}$ agrees with a bounded-Lipschitz metric dual to $L$.
\end{proposition}

\begin{proof}
Since those results are known when $A$ is unital, we assume that $A$ does not
have a unit. Let us assume that the Kantorovich distance $\kappa_{L}$ gives
$S(A)$ a finite diameter, say $r$. Let $a\in A$ be selfadjoint such that
$L(a)\leq1$. We wish to prove that $\left\Vert a\right\Vert \leq r$. Since $A$
is not unital, the weak*\ closure of $S(A)$ is $Q(A)$ and in particular it
contains the zero functional. Let $\varepsilon>0$. We can find $\varphi
_{\varepsilon}$ such that $\left\vert \varphi_{\varepsilon}(a)\right\vert
\leq\varepsilon$. Now, since $a$ is self-adjoint:%
\[
\left\Vert a\right\Vert =\sup_{\varphi\in S(A)}\left\vert \varphi
(a)\right\vert \leq\varepsilon+\sup_{\varphi\in S(A)}\left\vert \varphi
(a)-\varphi_{\varepsilon}(a)\right\vert \leq\varepsilon+\sup_{\varphi\in
S(A)}\kappa_{L}(\varphi,\varphi_{\varepsilon})
\]
since $L(a)\leq1$. Hence for all $\varepsilon>0$ we have proven $\left\Vert
a\right\Vert \leq r+\varepsilon$ so $\left\Vert a\right\Vert \leq r$ as we
wished. By definition, the distance $\kappa_{L}$ agrees in this case with
$d_{L,r,1}$.
\end{proof}

We conclude with some examples to illustrate Theorem (\ref{quasiLipalg}). Our
family of examples is rather simple, but already illustrates that our results
are natural:

\begin{example}
\label{compact}Let $\left(  X,\rho\right)  $ be a locally compact separable
metric space and $\mathcal{K}$ be the C*-algebra of compact operators. We let
$\mathcal{A}_{0}=C_{0}(X,\mathcal{K})$ be the C*-algebra of continuous
functions from $X$ into $\mathcal{K}$ and vanishing at infinity. The most
immediate choice for a Lipschitz seminorm $l$ on $\mathcal{A}_{0}$ is defined
by:%
\[
\forall a\in\mathcal{A}_{0}\ \ \ l(a)=\sup_{x,y\in X\text{,}x\not =y}%
\frac{\left\Vert a(x)-a(y)\right\Vert }{\rho(x,y)}\text{.}%
\]
Hence, we set $L=\max\left\{  l,\left\Vert .\right\Vert \right\}  $ as the
natural choice for our quasi-Lip-norm on $\mathcal{A}_{0}$, and indeed:

\begin{proposition}
Let $\left(  X,\rho\right)  $ be a locally compact separable metric space, and
let $l$ be the Lipschitz seminorm defined by $\rho$ on $C_{0}(X)$. The
seminorm $L=\max\left\{  l,\left\Vert .\right\Vert \right\}  $ is a
quasi-Lip-norm on $C_{0}(X)$, where $\left\Vert .\right\Vert $ is the supremum
norm on $C_{0}(X)$.
\end{proposition}

\begin{proof}
We see $\mathcal{K}$ as acting on the separable Hilbert space $\mathcal{H}$
with inner product $\left\langle .,.\right\rangle $ and we choose a Hilbert
basis $\left(  e_{n}\right)  _{n\in\mathbb{N}}$. For each $n\in\mathbb{N}$ we
set $P_{n}$ to be the orthogonal projection on the linear span of $\left\{
e_{0},\ldots,e_{n}\right\}  $. Of course, $\left(  P_{n}\right)
_{n\in\mathbb{N}}$ satisfies Condition (\ref{inductiveunits}). Moreover
$P_{n}\mathcal{A}_{0}P_{n}=C_{0}(X,M_{n+1})$ where $M_{n}$ is the C*-algebra
of $n\times n$-complex matrices. Let $f\in\mathcal{A}_{0}$ such that
$L(f)\leq1$ and $f=f^{\ast}$. Then, for all $k,l\in\mathbb{N}$:%
\[
1\geq\sup\left\{  \frac{\left\langle f(x)e_{k}-f(y)e_{k},e_{l}\right\rangle
}{\rho(x,y)}:x\not =y\right\}
\]
and $\left\Vert \left\langle fe_{k},e_{l}\right\rangle \right\Vert _{C_{0}%
(X)}\leq1$, so $\left\langle fe_{k},e_{l}\right\rangle $ is a 1-Lipschitz map
(for the distance $\rho$) valued in $[-1,1]$. By Arz\'{e}la-Ascoli, this
implies that $\left\langle \mathcal{B}_{L}e_{k},e_{l}\right\rangle $ is norm
totally bounded in $C_{0}(X)$. Hence $P_{n}\mathcal{B}_{L}P_{n}$ is norm
totally bounded in $C_{0}(X,M_{n+1})$ (by a diagonal argument). Therefore
$\mathcal{B}_{L}$ is \textrm{wu}-totally bounded.
\end{proof}

The construction of $L$ reflects the idea that $C_{0}(X)$ and $\mathcal{A}%
_{0}$ should have the same geometry (Morita equivalence). In particular, if
$\mathbb{Z}$ acts properly and freely on a locally compact space $Y$ so that
the orbit space is $X$ then $\mathbb{Z}\ltimes C_{0}(Y)$ is Morita equivalent
to $\mathcal{A}_{0}$. Moreover, if $H_{3}(X,\mathbb{Z})=0\ $then
$\mathbb{Z}\ltimes C_{0}(Y)$ is $\mathcal{A}_{0}$ by \cite[Corollary
10.9.6]{Dixmier}. If $d$ is a metric on $Y$ then one can define $\rho$ on $X$
by $\rho\left(  \mathcal{O}(x),\mathcal{O}(y)\right)  =\inf\left\{  d(z\cdot
x,z^{\prime}\cdot y):z,z^{\prime}\in\mathbb{Z}\right\}  $ (where
$\mathcal{O}(x)$ is the orbit of $x\in Y$), so that such crossed-products fit
in this example.
\end{example}

Example (\ref{compact}), though rather simple, contains in fact an important
lesson: indeed, let us use the notations of the proof of Proposition
(\ref{ctrexmp}). The sequence $\left(  S_{n}^{\ast}\right)  _{n\in\mathbb{N}}$
converges to $0$ for \textrm{wu}, though it is immediate that $\left(
S_{n}^{\ast}\right)  _{n\in\mathbb{N}}$ does not have any subsequence which
converges (to $0$) in the strict topology. First, simply by setting
$X=\left\{  1\right\}  $ in Example\ (\ref{compact}), we observe that
$L(S_{n})=\left\Vert S_{n}\right\Vert =1$ and this is counterexample to the
trivial fact that $\left\{  a\in\mathcal{K}:\left\Vert a\right\Vert
\leq1\right\}  $ is not strictly totally bounded, though it is \textrm{wu}
totally bounded. On the other hand, if $(X,\rho)$ is an arbitrary separable
locally compact metric space, we can pick $f\in C_{0}(X)$ such that
$\left\Vert f\right\Vert \leq1$ and which is $1$-Lipschitz for $\rho$. We then
define $a_{n}\in C_{0}(X,\mathcal{K})$ by $a_{n}=f\otimes S_{n}^{\ast}$. It is
easy to check that $L(a_{n})\leq1$ and $\left\Vert a_{n}\right\Vert \leq1$ for
all $n\in\mathbb{N}$, yet again $\left(  a_{n}\right)  _{n\in\mathbb{N}}$ has
not strictly convergent subsequence. Thus, $\left\{  a\in\mathcal{A}%
:L(a)\leq1\text{ and }\left\Vert a\right\Vert \leq1\right\}  $ is
\textrm{wu}-totally bounded yet not strictly totally bounded, even for such a
natural construction of a Lipschitz norm.

\bibliographystyle{amsplain}
\bibliography{thesis}

\end{document}